%

\magnification=1200
\input amstex
\documentstyle{amsppt}

\newcount\PLv\newcount\PLw\newcount\PLx\newcount\PLy\newdimen\PLyy\newdimen\PLX
\newbox\PLdot \setbox\PLdot\hbox{.} \def\scl{.08} 
\def\PLot#1{\PLx`#1\advance\PLx-42\PLy\PLx\PLv\PLx\divide\PLy9\PLw\PLy\multiply
\PLw9\advance\PLx-\PLw\advance\PLx-4\PLy-\PLy\advance\PLy4\PLX=\the\PLx
pt
\advance\PLyy\the\PLy pt\wd\PLdot=\scl\PLX\raise\scl\PLyy\copy\PLdot}
\def\draw#1{\ifx#1\end\let\next=\relax\else\PLot#1\let\next=\draw\fi\next}

\def\invamp{\hbox{\PLyy=70pt\draw
:::;DMV_gqppyyyyyooooxxxnnwvlutkjaWNE=5-./9%
9:::CCCC:::99/..--544=EENWWaajjjkktttttttNNNVVVVVVVV\end \hskip4pt}}
\newbox\iabox\setbox\iabox\invamp \def\Invamp{\copy\iabox}
\define\pa{\Invamp\,}
\define\ti{\otimes}
\define\li{\multimap}
\topmatter
\title
Is Game Semantics Necessary?
\endtitle
\author
Andreas Blass
\endauthor
\address
Mathematics Dept., University of Michigan, Ann Arbor, 
MI 48109, U.S.A.
\endaddress
\email
ablass\@umich.edu
\endemail
\thanks Partially supported by NSF grant DMS-9204276.
\endthanks
\abstract
We discuss the extent to which game semantics is implicit in the basic
concepts of linear logic.
\endabstract
\endtopmatter
\document

\head
Introduction
\endhead

The purpose of this paper is to show that a version of game semantics
for linear logic is implicit in the logic itself and the basic
intuitions underlying the logic.  Like the talk at CSL'93 on which it
is based, the body of this paper is intended to be accessible to
people with little or no previous knowledge of linear logic or game
semantics.  Comments that do presuppose such prior knowledge have been
relegated to a series of notes at the end of the paper.

\head
Propositions as Types
\endhead

The relevance of various constructive propositional logics, including
linear logic, to computation and particularly to type theory is
largely based on the propositions-as-types paradigm, also often called
the Curry-Howard isomorphism \cite{8, 9, 13}.  In its simplest form, this
paradigm involves a correspondence between the constructive logic of
implication and simple typed combinatory logic.  Constructive logic of
implication can be axiomatized by the schemes
$$
A\to(B\to A)
$$
and
$$
[A\to(B\to C)]\,\to\,[(A\to B)\to(A\to C)]
$$
and the rule of modus ponens
$$
\frac{A\to B\qquad A}{B}.
$$
Curry and Howard \cite9 noticed that, if one reads the letters as referring to
sets (or types) rather than propositions and reads $A\to B$ as the set
of functions from $A$ to $B$ rather than implication, then the two 
schemes are the types of the basic combinators $K$ and $S$, defined by 
$$
(Kx)y=x\quad\text{and}\quad ((Sx)y)z=(xz)(yz),
$$
and the rule of modus ponens corresponds to the typing of the basic
construction of combinatory logic, application of functions to
arguments.  All combinators definable in the simple typed
lambda-calculus can be obtained from $K$ and $S$ by repeated
application; this is the type-theoretic analog of the fact that all
constructively valid formulas involving only implication can be
obtained from the two schemes above by repeated use of modus ponens.

The correspondence between propositions and types can easily be
extended to include other connectives on the propositional side and
more constructions than the simple typed lambda-calculus on the types
side.  In particular, conjunction and disjunction of propositions
correspond to the cartesian product and the disjoint union of types,
respectively.  

Very similar ideas are contained in the intended interpretation of the
propositional connectives in intuitionistic mathematics \cite{4, 8, 14}.
There, the meaning of a proposition is specified by telling what is
required in order to prove the proposition, and connectives are
explained by telling how they affect proofs.  Specifically, a proof of
$A\&B$ is an ordered pair consisting of a proof of $A$ and a proof of
$B$, a proof of $A\vee B$ is a proof of $A$ or a proof of $B$ together
with the information which disjunct is being proved, and a proof of
$A\to B$ is a construction that converts any proof of $A$ into a proof
of $B$.  If one identifies a proposition with the set of its proofs,
then this Brouwer-Heyting explanation of the connectives exactly
matches the propositions-as-types interpretation described above, at
least if we regard ``construction that converts'' as an
intuitionistic way of referring to functions.

\head
Linearity
\endhead

One of the (two) fundamental ingredients of Girard's linear logic
\cite6 is a computational refinement of the type-theoretic notion of
function, paying attention to how many times the input is used in
computing the output.  In both classical and constructive mathematics,
the argument of a function is regarded as being permanently available
for arbitrarily repeated access during the computation of the
function's value at that argument, but it is clear that, if one is
interested in the efficiency of computation, it can be useful to know
how often an argument needs to be accessed (and even whether it is
accessed at all).  Motivated by these (and other) considerations,
Girard introduced $A\li B$ as the type of functions from $A$ to $B$ in
whose computation the argument in $A$ is accessed exactly once; such
functions are called {\sl linear\/}.  (This
sort of access counting is hard to define if one thinks of functions
in classical terms, but Girard showed that it makes good sense in
suitable computational situations, specifically for a special kind of
Scott domain called a coherence space.)  Girard also associated to
each type $A$ another type $!A$, the type of arbitrarily accessible
objects (or permanently stored objects) of type $A$.  A single access
of $!A$ consists of an arbitrary number of accesses of $A$.  Thus, the
traditional function space $A\to B$ can be described in this linear
framework as $(!A)\li B$.

Girard developed a logical system \cite6, a sequent calculus for
linear logic, which can be regarded as a variation of the standard
sequent calculi for propositional logic but which becomes much more
intuitive if viewed as being about types in the sense just explained
(and to be explained further below) rather than about propositions.

When one pays attention, as in linear logic, to the number of times a
data object is accessed, one discovers an ambiguity in the explanation
above of conjunction as cartesian product.  A data element of type
$A\&B$ should consist of a data element of type $A$ and one of type
$B$.  But does one access to $A\&B$ yield both components of the
ordered pair or just one?  At first sight, ``both'' might seem the
more reasonable answer, and in many contexts it is, but there are also
situations where one needs a type such that, in accessing it, one
specifies one of $A$ and $B$ and receives a data element of the
specified type.  For example, it is this sort of conjunction, not the
``both'' version, that yields a product in a category of types and
linear functions. Therefore, Girard included both versions of
conjunction (or product) in his system, using the notation $A\ti B$
for the ``both'' version and $A\&B$ for the ``only one'' version.
(See Note 1.  The former connective is usually called ``times'' 
and the latter ``with''.) The
logical rules of inference for these two connectives are
$$
\frac{\Gamma\vdash A\qquad\Delta\vdash B}{\Gamma,\Delta\vdash A\ti B}
\qquad\text{and}\qquad
\frac{\Gamma\vdash A\qquad\Gamma\vdash B}{\Gamma\vdash A\&B},
$$
where $\Gamma$ and $\Delta$ are finite lists of hypotheses.
Both of these are commonly used as introduction rules for conjunction
in sequent calculi, and they are equivalent in the presence of the
structural rules of thinning and contraction.  But these two
structural rules are not admissible in linear logic, for their effect
is precisely to alter the number of occurrences of a hypothesis, a number
to which linear logic insists on paying attention, as it corresponds
to the number of times input data are accessed.  Thus, in linear
logic, $\ti$ and $\&$ are genuinely different connectives.

\head
Access Protocols
\endhead

The introduction of the \& operation on types forces a revision in the
basic intuition of what a data type is.  Until now, we have regarded
types as essentially synonymous with sets.  Computationally, a type
could be regarded as a server from which a client can get, in one
access, an element of that type; the client need not do anything more
than show up.  We shall call such types {\sl simple\/} to distinguish
them from the more general types about to be introduced. 
(I thank Dexter Kozen for suggesting the
``client-server'' terminology; in my lecture I had talked about a
``user'' and a ``data resource''.)  For a type of the form $A\&B$,
however, the situation is different, for the client must specify
whether he wants an element from $A$ or one from $B$.  (For the
disjunction which, following Girard, we write $A\oplus B$, the client
need not do anything, for the server will choose one of $A$ and $B$,
announce which it chose, and then provide a data element.)  Thus, \&
types provide a first example of a non-trivial {\sl access
protocol\/} and a non-simple type.  

Calling the transmission of one bit (the choice of $A$ or $B$) a
protocol may seem like undesirable jargon, but once the client has to
do anything at all we soon arrive at situations where the word
``protocol'' seems quite reasonable.  Consider, for example, what
happens when a client accesses a server for the data type
$$
[(A\&B)\oplus(C\&D)]\&[(E\&F)\oplus(G\&H)].
$$
First, the client must specify which of the two types in square
brackets he wishes to access, as the overall data type is constructed
from these by \&.  Suppose the client chooses to access the second
component, $(E\&F)\oplus(G\&H)$.  Then, as this is a disjunction
(disjoint union), the server will choose one of the disjuncts, announce
its choice, and provide an element of that disjunct.  Suppose it
chooses $E\&F$.  After announcing that choice but before providing a
data element, the server must await the client's decision whether to
access $E$ or $F$, as these are combined with $\&$.  Thus the access
protocol in this case consists of three bits transmitted between
client and server before any actual data are provided.  One can
obviously build similar examples with longer access protocols.

In view of this situation, we regard a data type as consisting not
only of a set of possible data elements but also as the access
protocol that is to be run before a data element is provided. Of
course simple data types, where the client merely shows up
and the server provides a data element, are included in this scheme;
they are types with trivial access protocols.  (This
generalization of the notion of data type seems to be quite
independent of the generalization by admitting partially defined data
elements that is the intuitive basis for domain theory.  See \cite{11}
for a combination of the two generalizations.)

We could even incorporate the transmission of data at the end of the
client-server interaction into the access protocol.  Instead of having
the server provide, after the access protocol is complete, an element
of a set $S$ of possible data, we can regard that last transmission
from the server as just another step of the access protocol, a choice
of one disjunct in an $S$-indexed disjunction.

Formally, a protocol is a pair consisting of a (1) non-empty set $H$
of finite sequences, the possible histories of the protocol up to any
point in its execution, such that every initial segment of a sequence
in $H$ is also in $H$ and such that no infinite sequence has all its
finite initial segments in $H$, and (2) a function $N$ from $H$ into
$\{c,s,t\}$.  The intention is that, for any history $h\in H$,
$N(h)=c$ if the next step after $h$ in the protocol is to be an action
of the client, $N(h)=s$ if the next step is by the server, and
$N(h)=t$ if the protocol is terminated at $h$.  We require, in accord
with this intention, that a node $h$ with $N(h)=t$ cannot be a proper
initial segment of any other node in $H$; no history can go past a
terminal condition.  (We do not require the converse.  It is
permissible for $N(h)$ to be $c$ or $s$ even if there are no proper
extensions of $h$ in $H$.  This situation would mean that the client
or server is expected to do something but cannot.  The simplest
example is the simple data type corresponding to the empty set, where
the server is expected to provide an element but cannot.)  Notice
that, by requiring that no infinite sequence have all initial segments
in $H$, we have chosen to consider only protocols that always end
after finitely many steps.  The theory could be expanded to allow or
even require infinite runs (cf. \cite2), but nothing of this sort
seems to be implicit in the formalism or the underlying intuitions of
linear logic. (But see the discussion of $!A$ below.)

As suggested by the title of this paper, the protocols considered here
can be viewed as games (or debates or dialogs) between the client and
the server \cite{1, 2, 10, 11}.  In this connection, the server is
usually called the proponent or player, and the client is called the
opponent.  The protocol specifies who is to move (see Note 2) and what
moves are legal at any point during a play of the game.  Our
protocols, unlike some versions of games \cite{1, 2} but like the
versions in \cite{3,11}, do not specify winners and losers, but it
seems reasonable to regard a server as ``winning'' if it succeeds in
running the entire protocol (including the final step of delivering
data) without ever being in a situation where it is expected to act
but cannot.

Having extended the notion of data type to include access protocols,
we must explain how the connectives of linear logic are to act on
these data types.  The so-called additive connectives, \& and
$\oplus$, are fairly easy to handle, as they were involved in
the introduction of access protocols.  Specifically, the data type
$A\&B$ can be described by saying that its access protocol begins with
a choice by the client of either $A$ or $B$ and that the rest of the
protocol (and the final data delivery, if that is construed separately
from the protocol) is to be exactly as in $A$ or in $B$, according to
the client's initial choice.  The description of $A\oplus B$ is exactly
the same, except that the initial choice of $A$ or $B$ is made by the
server rather than by the client.

The description of $A\ti B$ is more complicated.  The basic intuition
is that the client is to get data of both types $A$ and $B$.  So the
client and server must carry out the access protocols for both $A$ and
$B$.  But should they be executed in parallel, or in sequence, or
interleaved? If parallel, then synchronously or asynchronously? If
interleaved then in what order?  The formalism of linear logic itself
provides some information about these questions.  For example, it is
provable in linear logic that 
$$
A\ti(B_1\oplus B_2)\,\vdash\,\dashv\,(A\ti B_1)\oplus(A\ti B_2).
$$
Here $\vdash\,\dashv$ means that sequents in both directions are
provable, and it follows that the left and right sides of such a
double sequent are equivalent in the deductive system; any formula
containing an instance of the left side is interdeducible with the
result of substituting the corresponding instance of the right side.
The left side of the displayed equivalence describes a $\ti$
combination in which one of the two components, $B_1\oplus B_2$,
begins with a choice by the server.  The equivalence says that this
combination amounts to a protocol in which the server begins by 
choosing either $A\ti B_1$ or $A\ti B_2$ and then the protocol for the
chosen constituent is executed.  That is, the first thing that should
happen in the protocol for $A\ti(B_1\oplus B_2)$ is a choice by the
server of $B_1$ or $B_2$ to replace $B_1\oplus B_2$ in the original
combination.  Summarizing: If one of the constituents of a $\ti$
combination has a protocol beginning with a choice by the server, then
the combination's protocol begins with the corresponding choice by the
server.  

What if both sides of $\ti$ have protocols beginning with choices by
the server?  The distributive law of linear logic displayed above
easily yields the equivalence
$$
(A_1\oplus A_2)\ti(B_1\oplus B_2)\,\vdash\,\dashv\,
(A_1\ti B_1)\oplus(A_1\ti B_2)\oplus(A_2\ti B_1)\oplus(A_2\ti B_2).
$$
This means that, if both protocols begin with choices by the server,
then, in the combined protocol, the server should begin by making
these choices in both components.  (See Note 3.)

Thus, in a $\ti$ combination, the
client need not do anything until it is his turn to act in both
constituents. What happens then is, however, not determined by the
formal system of linear logic.  There seem to be two reasonable
options.
\roster
\item The client chooses one of the two constituents and performs the
first action required by its protocol.
\item The client begins the protocols of both constituents.
\endroster
In either case, after the client performs the first action it will (in
general) be the server's turn to act in one or both protocols, and the
client would wait for the server's action(s) before proceeding, as
indicated by the distributive laws in the preceding paragraphs.  Also,
in option \therosteritem1, it is to be understood that the client is
obligated to make his choices, at each step where he is to act in both
components, in such a way that the protocols for both components are
ultimately completed. (See Note 4.)
 Another way to say this is that the protocol
for $A\ti B$ is terminated only when both sub-protocols are.  If one
sub-protocol is terminated and the other is not, then client and
server continue to run the latter protocol.

There are several reasons for preferring option \therosteritem1, but
they do not seem entirely conclusive.  One is that it provides greater
flexibility to the client, who can act in one protocol and wait for
the server's response there before committing himself to a particular
action in the other protocol.  Regarding $\ti$ as a very generous way
(compared to \&) to supply data, we view such flexibility as natural.

Perhaps a stronger argument concerns the ``translations'' of the two
options into the formalism of linear logic, analogous to the
distributive laws above.  These formalizations read
$$
(A_1\&A_2)\ti(B_1\&B_2)\,\vdash\,\dashv\,X
$$
where for option \therosteritem1 $X$ is
$$
X_1=[(A_1\&A_2)\ti B_1]\&[(A_1\&A_2)\ti B_2]\&
[A_1\ti(B_1\&B_2)]\&[A_1\ti(B_1\&B_2)]
$$
while for option \therosteritem2 $X$ is
$$
X_2=(A_1\ti B_1)\&(A_2\ti B_1)\&(A_1\ti B_2)\&(A_2\ti B_2).
$$
Neither of these equivalences is provable in linear logic, though the
direction $\vdash$ is provable in both cases.  In fact, we have
$(A_1\&A_2)\ti(B_1\&B_2)\,\vdash X_1$ and $X_1\vdash X_2$, but neither
implication is reversible.  To adopt either of the options above
therefore involves going beyond what is provable in linear logic, as
it involves adding an implication from $X_1$ or $X_2$ to
$(A_1\&A_2)\ti(B_1\&B_2)$. Since $X_1\vdash X_2$, option
\therosteritem1 involves less of an addition to linear logic than
option \therosteritem2.  In fact, if we were to adopt the implication
corresponding to option \therosteritem2, then the implication
corresponding to \therosteritem1 would be deducible.  

Finally, we might give an argument ``from consensus'' for option
\therosteritem1, namely that this option was adopted by all authors on
game semantics \cite{1, 2, 3, 10, 11}.  (Additional arguments will
arise when we consider duality in the next section.)

We turn next to the description of the data type $A\li B$.  Recall
that this is the type of linear functions from $A$ to $B$.  Thus the
following description of its access protocol (and data delivery) seems
reasonable.  A server for $A\li B$ acts like a server of type $B$
(intuitively, it supplies a value of type $B$), provided it has access
(once, as a client) to a server of type $A$ (intuitively, it is given
a value of type $A$).  The client of $A\li B$ must not only act as a
client of type $B$ (to whom the server will ultimately supply the
desired value) but must also provide (or act as) a server of type $A$,
supplying the data from $A$ required as input.  Thus, the client and
server of $A\li B$ play the same roles in the sub-protocol for $B$ but
opposite roles in the sub-protocol for $A$.  Furthermore, the server
should be permitted to switch from one sub-protocol to another
whenever it wishes; this amounts to allowing a function-evaluator
access to its argument (only once but) at a time chosen by the
evaluator.  It is fairly clear that denying the server such permission
would unduly restrict it even in simple situations.

Finally, $!A$ should be a re-usable version of $A$.  One access to
$!A$ should be an arbitrary number of accesses to $A$, the arbitrary
number being chosen by the client.  Thus, a run of the access protocol
for $!A$ should consist of several runs of the access protocol of $A$,
the initiation of new sub-runs being at the discretion of the client.
It is reasonable to require that, as long as the client makes the same
choices in two of the sub-runs, so does the server.  This requirement
corresponds to the intuition that a data element of type $!A$ is a
single, stored element of type $A$, so that repeated accesses give the
same result.  The formalism of linear logic does not demand this sort
of consistency --- it would permit a form of storage resulting in a
different element at each access --- but we regard this as contrary to
the intuitions underlying the $!$ concept. (See Note 5.)  We also
allow the client to switch freely from one sub-run to another and to
resume previously abandoned sub-runs.

It is not clear when a run of the protocol $!A$ should be regarded as
terminated, since the client could always start a new run of the
sub-protocol $A$.  A sensible convention, in view of the idea that the
client of $!A$ can access $A$ as often as he wishes, is that a run of
$!A$ is terminated only when the client declares it to be terminated.
This convention, unfortunately, allows infinite sequences in which
every finite initial segment is a permissible history of $!A$, with
the client never declaring termination.  There are several ways around
this difficulty, none of which I consider entirely satisfactory in the
present context of extracting a semantics from the formalism and
underlying intuitions of linear logic.  One possibility is to simply
allow protocols to have infinite runs.  Another is to work with games,
where there are winners and losers, rather than protocols, and to
declare any play of $!A$  which the client never terminates as won by
the server (as a penalty for the client's cheating).  And a third is
to require the client to announce, at the beginning of the protocol
for $!A$, how many accesses of $A$ he will use; then the $!A$ protocol
would be terminated when the announced number of $A$ sub-protocols
have been terminated.

With the interpretations above for $\li$ and $!$, we obtain the
following, quite reasonable interpretation of $A\to B$, which we
recall was expressed in linear logic as $(!A)\li B$.  A server of type
$A\to B$ is prepared to act as a server of type $B$ provided it has
repeated access (as client) to a consistent server of type $A$.

\head
Negation and Duality
\endhead

The second fundamental idea of linear logic (after the access counting
represented by linearity) is linear negation or duality.  The
intuition behind linear negation (and other connectives derived from
it) is considerably less clear than that behind the connectives
discussed in the preceding sections.  When discussing linear logic in
terms of a flow of questions and answers in a network, Girard \cite{6, 7}
has indicated that questions and answers of type $A$ are to be
regarded as answers and questions, respectively, of the negated type,
written $A^\perp$.  It seems natural to identify the questions and
answers in Girard's discussion with the actions of the client and the
server, respectively, in our protocols.  Thus, linear negation amounts
to interchanging the roles of the client and the server in a protocol.
In this connection, it is convenient to adopt the viewpoint, mentioned
earlier, that the actual delivery of data by the server at the end of
the protocol is considered as just another piece of the protocol.
Thus, in the negated data type, this final action would be performed
by the client.

It is part of the formalism of linear logic \cite6 that negation is
involutive, i.e., $A^{\perp\perp}=A$.  In addition, Girard introduced
De~Morgan-style duals $\pa$ and $?$ for $\ti$ and $!$, respectively,
i.e., 
$$
A\pa B=(A^\perp\ti B^\perp)^\perp\qquad ?A=(!(A^\perp))^\perp.
$$
The connectives \& and $\oplus$ are dual to each other in the same
sense, so
$$
(A\oplus B)^\perp=A^\perp\&B^\perp\qquad
(A\&B)^\perp=A^\perp\oplus B^\perp.
$$
Also, in Girard's system, the linear implication can be expressed in
terms of $\ti$ and negation as
$$
A\li B=(A\ti(B^\perp))^\perp=A^\perp\pa B.
$$
To the extent that they involve only negation and the connectives
discussed in previous sections, these equations are correct when read
as descriptions of protocols.  For example, $(A\oplus
B)^\perp=A^\perp\&B^\perp$ amounts to the fact that $\oplus$ is
interpreted just like $\&$ except that the first choice is made by the
server instead of the client.  The correctness of the protocol reading
of $A\li B=(A\ti(B^\perp))^\perp$ constitutes additional support for
our decision to use option \therosteritem1 in designing the protocol
for $\ti$.  Had we used option \therosteritem2, then to keep $A\li
B=(A\ti(B^\perp))^\perp$ correct we would have to modify the protocol
for $\li$ by imposing an unnatural synchronization requirement on how
the server runs the two sub-protocols.

The new connectives $\pa$ and $?$ defined above seem rather unnatural
in computational terms, as they allow the server to make
choices that servers don't ordinarily make.  For example, the protocol
for $A\pa B$ consists of interleaved runs of the protocols for $A$ and
for $B$, with the client required to act in one or both protocols
whenever he can and the server required to act in one protocol
whenever it is due to act in both.  Thus, whenever the server acts, it
can choose which sub-protocol to act in (unless one is already
terminated).  Since our definition of $\ti$ required the client to
ultimately finish both subprotocols, the same requirement is
automatically imposed on the server in a $\pa$ combination.  In other
words, the server in $A\pa B$ will ultimately have to provide data of
both types $A$ and $B$.  Thus,
although it is De~Morgan dual to a sort of conjunction, $\ti$, the
connective $\pa$ is computationally more like a conjunction than like
a disjunction.

By duality, the distributive law for $\ti$ over $\oplus$ implies a
distributive law for $\pa$ over $\&$.  This fact leads to another (weak)
argument against option \therosteritem2 in the interpretation of
$\ti$.  Indeed, the sequent formulation of that option, namely
$(A_1\&A_2)\ti(B_1\&B_2)\,\vdash\,\dashv\,X_2$, says that $\ti$ behaves,
with regard to a particular instance of distributivity, as if it were
an entirely different connective $\pa$.  (Incidentally, these
distributive laws were the reason for Girard's notation for the
connectives, in which dual connectives do not have analogous symbols.
They also provide one motivation for the terminology ``additive'' for
the connectives \& and $\oplus$ and ``multiplicative'' for $\ti$ and
$\pa$, but I believe the original motivation for this came from the
coherence space interpretation of these connectives.) 

We close this section by recording for reference the nullary analogs
of the binary connectives $\&$, $\oplus$, $\ti$, and $\pa$.  Girard's
notation for these constants is $\top$, 0, 1, and $\bot$ respectively.
The protocol for $\top$ is that the client is expected to act first
but no action is possible.  (In any \& combination, the client acts
first and chooses one of the types being combined by \&; our
description of $\top$ is the special case where the number of types
being combined is zero.)  So $\top$ is a data type that simply cannot
be accessed.  Dually, 0 is a data type in which the server is expected
to act first but has no possible action; it is the empty (simple,
i.e., with trivial access protocol) data type.  The other two
constants, 1 and $\bot$, are both interpreted as data types for which
the access protocol is vacuous, i.e., the protocol is already terminated
without either participant doing anything.  (This protocol is
therefore even shorter than the simple ones, for in
the latter the client does nothing but the server is expected to
provide a data element.)  If, as is customary in linear logic, one
regards $\pa$ as a sort of disjunction, then $\bot$ is a sort of
``false'' (the empty disjunction), while 1 is a sort of ``true'' (the
empty conjunction), so giving them the same interpretation seems very
unreasonable.  It becomes reasonable in interpretations oriented more
toward computations than toward logic (more toward types than
propositions); indeed this identification of 1 and $\bot$ occurs both in
Girard's coherence space semantics
\cite6 and in the Abramsky-Jagadeesan version of game semantics
\cite1.  It becomes even more reasonable in our present situation
since, as pointed out above, $\pa$, though logically like a
disjunction, is computationally more like a conjunction.

\head
Truth and Validity
\endhead

In the preceding sections, we have described how to interpret as a
protocol any combination of atomic formulas, built using the
connectives of linear logic, provided an interpretation of the atomic
subformulas as protocols is given.  The interpretation extends easily
to sequents, for in linear logic a one-sided sequent $\vdash
A_1,A_2,\dots,A_n$ is treated exactly like the formula
$A_1\pa A_2\pa\dots\pa A_n$, and a two-sided sequent
$A_1,\dots,A_m\vdash B_1,\dots,B_n$ is treated exactly like the
one-sided sequent $\vdash{A_1\!}^\perp,\dots,{A_m\!}^\perp,B_1,\dots,B_n$.

We have not yet said anything about validity of formulas or sequents.
(By the preceding equivalences, it suffices to consider only
formulas.)  It is as if, in many-valued logic, we had provided the
truth tables for the connectives but had not indicated which truth
values are the distinguished ones.  We have postponed the issue of
defining validity because linear logic seems to give us less guidance
here than in interpreting the connectives and because we find a
greater divergence between the logical and computational points of
view. 

From the point of view of logic, validity of a formula should mean
that, no matter how its atomic subformulas are interpreted as
protocols, the protocol denoted by the formula is true.  This merely
shifts the problem to defining ``true'' for protocols, but here the
Brouwer-Heyting description of intuitionistic truth can provide some
guidance.  That description calls a formula true if and only if there
is a proof of it.  Thinking of propositions as types, in this case
simple types of proofs, we find truth identified with non-emptyness of
a simple data type.  So in this special situation, $A$ is true if and
only if the server can run the protocol for $A$ without getting stuck
(i.e., without being expected to act and having no action available).
This description makes sense also for more elaborate protocols, so one
could identify truth of a protocol with the existence of a behavior
of the server that contains answers whenever required by the protocol.
In game-theoretic language, this amounts to a winning strategy for the
server in the game where the server's objective is just to execute the
protocol without getting stuck.  Formally, a {\sl behavior\/} for the
server in a protocol $(H,N)$ is a subset $B$ of the set $H$ of
histories such that 
\roster
\item The empty sequence is in $B$.
\item If $h\in B$ with $N(h)=c$, then every one-term extension of $h$
is also in $B$.
\item If $h\in B$ with $N(h)=s$, then exactly one one-term extension
of $h$ is in $B$.
\endroster
Intuitively, $B$ is the set of histories that can arise when the
server behaves in a particular way.
(We could also consider non-deterministic behaviors, or partial
strategies, where ``exactly one'' is replaced with ``at least one'' in
\therosteritem3.)  Thus, we regard a formula (or sequent) as true in a
particular interpretation if it admits a behavior for the server, and
we regard it as valid if it is true in every interpretation.  

This approach to validity, motivated by standard ideas from logic, was
used, for example, in \cite2.  It works well in the context of
infinite games (as in \cite2), but seems seriously deficient when
games are required to terminate after a finite number of moves.  The
reason is that, by a classical theorem of Gale and Stewart \cite5,
such finitely long games always admit winning strategies for one or
the other player.  This means that, for any protocol $A$, either there
is a behavior (a winning strategy in the sense mentioned above) for
the server in $A$ or there is a winning strategy for the client, i.e.,
a strategy by which the client can force the server to get stuck.  But
such a strategy for the client is a winning strategy for the server in
$A^\perp$.  Thus, for any $A$, the server has a behavior in the
protocol $A^\perp\oplus A$, namely to initially choose whichever of
$A$ and $A^\perp$ has a winning strategy for the server and then to
follow that strategy.  Thus $A^\perp\oplus A$ is valid in the sense of
the preceding paragraph.  I regard this as a deficiency of the
semantics because $A^\perp\oplus A$ seems quite unreasonable in the
context of linear logic, both at the intuitive level and at the formal
level.  To make the point more explicit, we point out that from
$A^\perp\oplus A$ and the rule of thinning one can deduce in linear
logic the rule of contraction (see Note 6), whereas these are normally
regarded as quite different matters, thinning being more innocuous
than contraction \cite7.

The behaviors that witness the validity of $A^\perp\oplus A$ are of a
rather complicated sort.  To determine its first move, the server must
completely analyze the game associated to $A$.  It seems reasonable to
require that validity entail the existence of simpler, more uniform
strategies, not depending on a detailed analysis of sub-games.
(Consider, for example, the behaviors for the provable formula
$A^\perp\pa A$, which consist of merely copying the client's actions
from one component to the other.)  Instead of asking that a formula
``(have a behavior) independently of the interpretation of atomic
subformulas'' we ask that it ``have (a behavior independent of the
interpretation of atomic subformulas).''  It is not entirely clear
what independence in the second sense should mean.  Abramsky and
Jagadeesan \cite1 have introduced a strong notion of independence,
namely that when additional moves are added to the atomic subgames,
the new strategy should be an extension of the old.  They show that
this requirement, in conjunction with a requirement of
history-freeness (meaning that the server's action at any stage should
depend only on the client's immediately previous action, not on the
earlier history) forces strategies to do nothing more complex than
immediately copying the client's moves between the various subgames.
Their work shows that such a restriction, which may seem too strong at
first sight, works
well in multiplicative linear logic.  They also observed, however, that
history-freeness does not work well in the presence of the additive
connectives.  One possible meaning for ``independence'' is the
Abramsky-Jagadeesan notion of uniformity (without history-freeness).
Another possibility, a modification of an idea of Lorenzen \cite{12}, is
that the server should never find itself expected to act in an atomic
subgame $A$ unless the client has previously had to act in a
corresponding $A^\perp$.  All these speculations need considerably
more work in the direction of soundness (and if possible completeness)
theorems.  I have not (yet) pursued this, partly because this line of
thought seems to be getting more into technical adjustments of the
semantics and farther from the theme of what the intuitions and
formalism of linear logic tell us about its semantics.

\head
Notes
\endhead

\subhead
Note 1
\endsubhead
In classical mechanics, a particle has position $\ti$ momentum.  In
quantum mechanics, a particle has position \& momentum.  (Of course,
this is an over-simplification, ignoring the quantitative trade-off,
given by the uncertainty principle, between partial information about
position and momentum in quantum mechanics.)

\subhead
Note 2
\endsubhead
Our protocols always specify which of the two participants is expected
to initiate the interaction.  In the terminology of \cite1, they have
definite polarity.  To obtain their completeness theorems for
multiplicative linear logic, Abramsky and Jagadeesan \cite1 and Hyland
and Ong \cite{10} made essential use of non-polar games, i.e., games that
either player can start (and that may have entirely different rules
depending on who starts).  As noted in \cite1, non-polarity does not
work well in the presence of the additive connectives, which explains
why we, starting from the behavior of additive connectives, arrived at
polar games.

\subhead
Note 3
\endsubhead
In their analysis of game semantics for multiplicative linear logic,
Abramsky and Jagadeesan \cite1 found that having the server begin both
components of $A\ti B$ when possible was the source of the excessive
supply of valid formulas (including thinning) in \cite2.  So in their
semantics, the server would act in only one component in such a
situation. The other component would be started (if at all) by the
client.  This set-up depends on having non-polar games, so that the
client can start the other component.  Thus, as indicated in Note 2,
it does not work as well in the presence of the additives.  Indeed,
the remarks about distributivity in the main text indicate that, if
one wishes to deal with $\ti$ and $\oplus$ together (and maintain the
correctness of linear deductions, in particular of the distributive
law), then one is forced to the convention that, when both sides of a
$\ti$ are to be started by the server then it should start them both
before expecting the client to do anything.

\subhead 
Note 4
\endsubhead
Requiring the client to finish both parts of $A\ti B$ will (once we
explain validity) make thinning invalid in this semantics of
protocols.  Thinning, which can be expressed as $A,B\vdash A$, is
valid in the semantics of \cite2 because the server (there called
proponent) can play a copying strategy between the $A$'s without ever
entering the subgame $B$.  The situation is different here, as the
server must enter the subgame $B$ and my find itself called upon to
move but having no move available there.

\subhead
Note 5
\endsubhead
Girard has pointed out \cite7 that, unlike the additive and
multiplicative connectives, the exponentials are not determined by
their introduction rules.  That is, one could introduce into the
sequent calculus a second (dual) pair of exponentials $!'$ and $?'$,
with the same introduction rules as the original ! and ?, and it would
not be provable that they are equivalent to the originals.  One
possibility for such multiple exponentials would be to interpret one
with a consistency requirement, as in the main text, and the other
without such a requirement.  Intuitively, the former describes storage
of a single data element which can be reliably retrieved arbitrarily
often, while the latter describes a source of a stream of data
elements which can be accessed repeatedly and may provide different
elements at each access.  Notice that the approximation of $!A$ by
iterated ``times'' of $1\&A$, suggested in \cite6, corresponds to the
stream picture, not to reliable storage.

\subhead
Note 6
\endsubhead
The deduction of $\vdash A,\Gamma$ from $\vdash A,A,\Gamma$ in the
presence of $\vdash A^\perp\oplus A$ and thinning proceeds as follows.
Obtain $\vdash A^\perp,A,\Gamma$ by thinning from the axiom $\vdash
A^\perp,A$.  Combine it with the assumption $\vdash A,A,\Gamma$ by the
\& rule to get $\vdash A\&A^\perp,A,\Gamma$.  Then cut the
$A\&A^\perp$ in this last sequent against the assumption
$\vdash A^\perp\oplus A$, to get the required $\vdash A,\Gamma$.

\enddocument